\documentclass[review]{elsarticle} 
            
\usepackage{amsmath}
\usepackage{graphicx}
\usepackage{booktabs}
\usepackage{enumitem}
\usepackage{amsfonts}
\usepackage{amssymb}
\usepackage{cancel}
\usepackage{mathtools}
\usepackage{float}
\usepackage{url}
\DeclareMathOperator*{\argmax}{arg\,max}

\usepackage{accents}
\usepackage{caption}
\usepackage{todonotes}
\usepackage{subcaption}
\usepackage{lineno,hyperref}
\PassOptionsToPackage{unicode}{hyperref}
\PassOptionsToPackage{naturalnames}{hyperref}
\usepackage{cleveref}

\begin{document}

\begin{frontmatter}

\title{Sensitivity Analysis of Passenger Behavioral Model for Dynamic Pricing of Shared Mobility on Demand}

\author[main]{Vineet Jagadeesan Nair\corref{mycorrespondingauthor}}
\cortext[mycorrespondingauthor]{Corresponding author}
\ead{jvineet9@mit.edu}

\author[main]{Yue Guan}
\ead{guany@mit.edu}

\author[main]{Anuradha M. Annaswamy}
\ead{aanna@mit.edu}

\author[ford]{H. Eric Tseng}
\ead{htseng@ford.com}

\author[ford]{Baljeet Singh}
\ead{bsing124@ford.com}

\address[main]{Department of Mechanical Engineering, Massachusetts Institute of Technology, Cambridge, MA 02139, USA}
\address[ford]{Research and Advanced Engineering (R\&A),
Ford Motor Company, 1 American Road, Dearborn, MI 48126, USA}

\begin{abstract} 
This paper provides a framework to quantify the sensitivity associated with behavioral models based on Cumulative Prospect Theory (CPT). These are used to design dynamic pricing strategies aimed at maximizing performance metrics of the Shared Mobility On Demand Service (SMoDS), as solutions to a constrained nonlinear optimization problem. We analyze the sensitivity of both the optimal tariff as well as the optimal objective function with respect to CPT model parameters. In addition to deriving analytical solutions under certain assumptions, more general numerical results are obtained via computational experiments and simulations to analyze the sensitivity. We find that the model is relatively robust for small to moderate parameter perturbations. Although some of the trends in sensitivity are fairly general, the exact nature of variations in many cases depends heavily on the specific travel scenarios and modes being considered. This is primarily due to the complex nonlinearities in the problem, as well as the significant heterogeneity in passenger preferences across different types of trips.
\end{abstract}


\end{frontmatter}


\section{INTRODUCTION}
\subsection{Motivation}
In recent times, several ride sharing platforms have emerged, that offer the potential for increased affordability, convenience and customizability \cite{ambrosino2004mobility}\cite{chong2013autonomy}. There is also a growing shift away from exclusive, door-to-door ridesharing services towards ride-pooling, which offers additional benefits including mitigating traffic congestion, reducing cumulative travel times and emissions while increasing fleet utilization rates \cite{santi2014quantifying}. This paper pertains to such Shared Mobility on Demand Services (SMoDS) i.e. large-scale pooled ridesharing. Our prior work described a comprehensive solution with dynamic routing via an alternating minimization (AltMin) optimization algorithm \cite{guan2019dynamic} and dynamic pricing using a passenger behavioral model based on Cumulative Prospect Theory (CPT) \cite{guan2019cumulative}. This paper is an extension exploring the sensitivity and robustness of the ridesharing system's performance to errors and uncertainty associated with the estimation and calibration of the CPT passenger behavioral model used to set dynamic tariffs. The model parameters are obtained by collecting data from users through a comprehensive survey involving discrete choice experiments \cite{ben1985discrete}, followed by maximum simulated likelihood estimation to obtain the population-level mode-choice model \cite{train2009discrete} and nonlinear least squares to arrive at individual specific CPT parameters describing their risk attitudes \cite{rieger2017estimating}. 


\subsection{Background and literature review}
\begin{figure}[H]
  \centering
    \includegraphics[width=\linewidth]{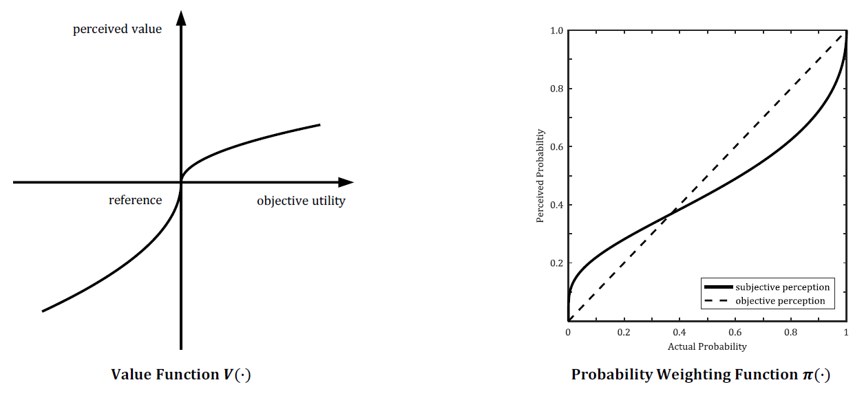}
  \caption{Illustrations of the CPT value function and probability weighting functions \cite{guan2019cumulative}.}
  \label{fig:CPTfuncs}
\end{figure}

Conventional Expected Utility Theory (EUT) postulates that consumers choose among travel options based on their respective expected utilities \cite{fishburn1988nonlinear}. Cumulative Prospect Theory (CPT) is an alternative to EUT that better describes subjective human decision making in the presence of uncertainty and risk \cite{tversky1992advances}\cite{kahneman2013prospect}. This is necessary in the case of SMoDS due to the significant variability in travel times introduced by pooling rides. The CPT behavior model used in this study is described by the value function $V(\cdot)$ and probability distortion $\pi(\cdot)$ given by \cite{guan2019cumulative}, with $\pi(0) = 0$ and $\pi(1) = 1$ by definition. These nonlinearities transform the objective utilities ($u$) and probabilities ($p$) of each possible outcome to subjective values, as perceived by the passengers. The graphs in \cref{fig:CPTfuncs} show examples of how the value and probability weighting functions may vary according to the objective utility $u$ and actual probability $p$, respectively.
\begin{align}
    V(u) & = \begin{cases} (u-R)^{\beta^+} &\mbox{if } u \geq R \\
    -\lambda(R-u)^{\beta^-} &\mbox{if } u < R
    \end{cases}
    \label{eq:value_func} \\
    \pi(p) & = e^{-(-ln(p))^\alpha} 
\end{align}

The CPT parameters here describe loss aversion ($\lambda$), diminishing sensitivity in gains ($\beta^+$) and losses ($\beta^\_$) and probability distortion ($\alpha$). The reference $R$ is the baseline against which users compare uncertain prospects. These can vary across individuals and also depending on the particular set of alternatives the customer is facing.

There are several sources of error and uncertainty involved in the CPT parameter estimation process that need to be addressed. These include sampling errors, survey design issues, response biases etc. Population level mode-choice models require very large sample sizes to obtain relatively accurate distributions for utility function coefficients and still involve some finite levels of uncertainty. Furthermore, travel preferences of consumers are fluid and can vary significantly over time, among different individuals in a population, as well as depending on the particular trip characteristics being considered and options in the current choice set. However, surveys capture only a static snapshot for a subset of passengers and a limited number of travel choice scenarios. Additionally, CPT parameters describing risk attitudes are often specific to each user and determined from a much smaller set of only their responses, rather than the whole sample of respondents. Thus, they lack statistical properties like asymptotic normality even with large sample sizes \cite{wang2019risk}. Inaccuracies in model parameters can result in setting sub-optimal dynamic tariffs that reduce the operational efficiency of the SMoDS, leading to decreases in revenue, ridership and fleet utilization rates. To summarize, the fact that these parameters are very passenger and situation-specific makes the model more prone to errors

Since the novel aspect of this behavioral model is the incorporation of CPT and because the parameters associated with it involve a greater deal of uncertainty compared to the mode-specific utility functions, this work will focus on analyzing the sensitivity only with respect to these CPT parameters. To our knowledge, CPT has not been explored for SMoDS in the literature and neither has such a sensitivity analysis been considered to date.

There has been significant past work conducted in analyzing the sensitivity of parametric linear and nonlinear optimization problems. The early work in this field analyzed the sensitivity and stability of nonlinear programs (NLP) \cite{fiacco1983introduction}\cite{fiacco1990sensitivity}. These established the mathematical foundation including basic theorems related to the smoothness, continuity and differentiability properties of the optimal solution and value function \cite{shapiro1988sensitivity}, performing both first and second order sensitivity analyses as well as determining asymptotic bounds on sensitivity derivatives \cite{shapiro1985second}. Several different approaches and solution techniques have also been explored, including penalty-function methods \cite{fiacco1976sensitivity}, generalized perturbation approaches \cite{castillo2006perturbation} and directional derivatives \cite{ralph1995directional}. In addition to localized analyses that focus on varying one parameter at a time, methods for global sensitivity analysis have also been studied, often using Monte Carlo techniques \cite{wagner1995global}. 

The main contribution of this paper is applying tools from sensitivity and robustness analysis specifically to behavioral models based on prospect theory, thereby developing tools to quantify the significant uncertainty and estimation errors associated with such models in the SMoDS context. This will enable us to develop corrective measures to mitigate the negative effects of such model uncertainty in future work.

\section{METHODOLOGY}
\subsection{\label{sec:prob} Problem formulation}
This study will focus on a single passenger taking a trip using the SMoDS. Without loss of generality, we consider a case where the passenger chooses between only two modes of transport: the uncertain SMoDS $sm$ against a baseline travel alternative $A$ (e.g. public transit, driving or exclusive door-to-door ridesharing) that can be treated as a certain prospect, offering a fixed objective utility $u_o$.
The objective utility of a trip with a certain travel option is calculated using a linear multinomial logit mode, based on the travel times spent on each leg, tariff $\gamma$ and the alternative-specific constant (ASC) $c$ of the service:
\begin{align}
    u & = \mathbf{a}^\intercal \;\mathbf{t} + b\gamma + c \\
    & = x + b\gamma
\end{align}
where $\mathbf{t} = [t_{walk}, t_{wait}, t_{ride}]^{\intercal}$ denote the walking, waiting, and riding times, respectively and $\mathbf{a}=[a_{walk}, a_{wait}, a_{ride}]^{\intercal}$ are the travel time coefficients for each leg. Here, $x$ is used to compactly represent the component of objective utility due to all travel times on different legs combined along with the ASC of the travel mode. The coefficients on travel times ($\mathbf{a}$) and tariff ($b$) are negative since these represent disutilities to the consumer, while $c$ can be either positive or negative depending on the characteristics of the given travel option. 

For simplicity and analytical tractability, for any given ride offer, the possible outcomes with the SMoDS are modelled as following a Bernoulli distribution, i.e., it is assumed to have only two possible travel time outcomes $\mathbf{\overline{t}}$ and $\mathbf{\underline{t}}$ ($\mathbf{\underline{t}} \leq \mathbf{\overline{t}}$) having corresponding utilities $\underline{u}$ and $\overline{u}$ ($\underline{u} \leq \overline{u}$), occurring with probabilities of $p \in [0,1]$ and $1-p$ respectively. This choice of distribution is a reasonable starting point and makes the problem more tractable. Moreover, accurately estimating probability distributions for travel times expected by each passenger would require a very large number of draws for each respondent in the population. Thus, any distribution fitted using data from a reasonably large sample size will still involve some finite error. We also assume that both outcomes offer the same price since most ridesharing services guarantee trip tariff at the time of ride offer. The analysis we present below can readily be extended to situations with multiple travel alternatives. This framework can also be applied to model with more complex probability distributions for SMoDS travel times, having more than two outcomes. These include both continuous (e.g. Gaussian, extreme value) and discrete (e.g. Poisson) distributions. The specific assumptions made here are primarily for simplicity and tractability while deriving analytical results, the broader insights and trends also hold for more general cases. The utilities of the SMoDS and alternative are then given by:
\begin{align}
\label{eq:utilities}
 \underline{u} & = \mathbf{a}_{sm}^\intercal \;\mathbf{\overline{t}} + b_{sm}\gamma_{sm} + c_{sm} = \underline{x} + b\gamma \\
 \overline{u} & = \mathbf{a}_{sm}^\intercal \;\mathbf{\underline{t}} + b_{sm}\gamma_{sm} + c_{sm} = \overline{x} + b\gamma\\
 u_o & = \mathbf{a}_o^\intercal \;\mathbf{t}_o + b_o\gamma_o + c_o
\end{align}
If $u_o \leq \underline{u}$, the customer would always choose the SMoDS since it offers strictly better outcomes and conversely if $u_o \geq \overline{u}$, they would always choose option $A$. Thus, the only cases considered are where $\underline{u} \leq u_o \leq \overline{u}$ are considered (note: $u_o$ can still be either a gain or loss) such that the consumer's choice (of accepting or rejecting the SMoDS ride offer) is non-trivial. Given that the SMoDS outcomes follow a Bernoulli distribution, its cumulative distribution function (CDF) is defined on the support $[\underline{u},\overline{u}]$:
\begin{equation}
    F_U(u) = \begin{cases} 0 &\mbox{if } u < \underline{u} \\
    p &\mbox{if } \underline{u} \leq u < \overline{u} \\
    1 &\mbox{if } u \geq \overline{u} \ \end{cases}
\end{equation}
\subsection{CPT model overview}
The study focuses on analyzing the model's sensitivity to the CPT model parameters as follows. The reader is referred to \cite{guan2019cumulative} for more details regarding the CPT based passenger behavioral model in the SMoDS context:
\begin{enumerate}
    \item \textbf{CPT parameters}: It assumed for simplicity that $\beta^+ = \beta^- = \beta$ i.e. the passenger displays similar reductions in sensitivity while moving away from the reference value, in both the gain and loss regime. The sensitivity with respect to these parameters is computed in the standard sense by allowing continuous variations in their values. 
    \item \textbf{Reference ($R$)}: Treated as a hyper-parameter and case studies using a few different reference types are considered: 
    \begin{itemize}
        \item \underline{Static}: Fixed for each customer, independent of the SMoDS ride offer. This could be set as the objective utility of the most frequently used travel alternative (excluding SMoDS) i.e. $R = u_o$.
        \item \underline{Dynamic}: $R$ depends on the uncertain prospect itself i.e. it varies with the SMoDS offer. Some examples considered here are:
        \begin{itemize}
            \item Expected utility of SMoDS 
            \begin{equation}
                R = p \underline{u} + (1-p) \overline{u}
            \end{equation}
            \item Best ($R = \overline{u}$) or worst-case ($R = \underline{u}$) utilities corresponding to the shortest ($\mathbf{\underline{t}}$) and longest ($\mathbf{\overline{t}}$) travel times, respectively. 
        \end{itemize}
    \end{itemize}
    \item \textbf{Probability distributions}: These are the probabilistic distributions of expected travel times spent on different legs (i.e. walking, waiting and riding), as perceived by the users themselves. In the current study, we model the objective utility of the SMoDS as a binary random variable. Thus, this distribution can be varied by altering the parameter $p$ which is the probability of the worst-case SMoDS outcome ($\underline{u}$) occurring.
\end{enumerate}
$U_R^s$ and $A^s_R$ are then the subjective utilities of the SMoDS and the alternative $A$, as perceived by the passenger relative to the reference point $R$. These can be calculated as follows:
\begin{align}
    A^s_R & = \pi(1) \cdot V(u_o)  = V(u_o) \\
    U^s_R & = w_1 \cdot V(\underline{u}) +  w_2 \cdot V(\overline{u})
\end{align}
where $w_1$ and $w_2$ are subjective probability weights calculated using:
\begin{align}
    w_i & = \begin{cases} \pi[F_U(u_i)] - \pi[F_U(u_{i-1})] &\mbox{} u_i < R \\
    \pi[1-F_U(u_{i-1})]- \pi[1-F_U(u_i)] &\mbox{}  u_i \; \geq R
    \end{cases}
    \label{eq:weights}
\end{align}
where $u_1 = \underline{u} \leq u_2 = \overline{u}$ and $F_U(\cdot)$ is the cumulative distribution function (CDF) of the SMoDS utility. The subjective probability of acceptance can be calculated using:
\begin{align}
    p^R_s(\gamma;\Vec{\theta}) & = \frac{e^{U^s_R}}{e^{U^s_R} + e^{A^s_R}} \label{eq:subj_prob}
\end{align}
\noindent where $\Vec{\theta} = [\alpha, \beta, \lambda, p, R]^{\intercal}$ consists of all the parameters of interest, assembled together. From eq. (\ref{eq:value_func})-(\ref{eq:subj_prob}), it is easy to see that $p_s^R$, the main output of the CPT-based passenger behavioral model is a function of the SMoDS tariff $\gamma$ parametrized by $\Vec{\theta}$. In the following, we will evaluate this model's sensitivity with respect to these key parameters by formulating the problem as one of nonlinear optimization.

\subsection{Optimization}
The dynamic tariff $\gamma$ is set by solving a constrained, nonlinear parametric optimization problem. We can consider several possible objective functions. For example, maximizing expected ridership for the fleet would be equivalent to directly maximizing acceptance probability $p^s_R$ itself. If we maximize expected revenue, the below NLP results:
\begin{align}
    \min_{\gamma} \quad & - f(\gamma; \Vec{\theta}) \triangleq - \gamma \cdot p_s^R(\gamma; \Vec{\theta}) \label{eq:exp_rev}\\
    \textrm{s.t. $g^1$:} \quad & \underline{\gamma} - \gamma \leq 0 \label{eq:opt_prob1} \\
    \textrm{$g^2$:} \quad & \gamma - \overline{\gamma}  \leq 0
    \label{eq:opt_prob2}
\end{align}
\noindent \cref{eq:exp_rev} reflects the expected revenue because the revenue per passenger per trip is $\gamma$ with a probability $p^s_R$ and $0$ with probability $1 - p^s_R$. Yet another option is to perform a weighted multiobjective optimization considering multiple objectives like revenue, ridership and utilization. The sensitivity analysis is then aimed at understanding how changes in parameters affect the optimal dynamic tariff $\gamma^*$ and value of the objective function $f(\gamma; \Vec{\theta})$, where $\gamma$ is the decision variable and $\Vec{\theta}$ represents all the model parameters. Note that the constraints in \cref{eq:opt_prob1} and \cref{eq:opt_prob2} only include upper and lower bound on the dynamic tariff charged such that it is within a reasonable range, i.e., $\gamma \in [\underline{\gamma}, \overline{\gamma}]$. All other constraints related to travel times have already been accounted for by the routing algorithm in generating the SMoDS ride offer and its possible outcomes. In practice, the lower bound could be the minimum break-even price per trip and the upper bound could be some sensible limit e.g. SMoDS tariff cannot be higher than that of exclusive ridesharing.
\subsection{Optimality conditions}
In the following, subscripts indicate partial derivatives w.r.t. that variable. The Lagrangian dual for the NLP formulated in eq. (\ref{eq:exp_rev})-(\ref{eq:opt_prob2}) is \cite{bertsekas1997nonlinear}:
\begin{equation}
    \mathcal{L}(\gamma; \Vec{\theta})  = -f(\gamma; \Vec{\theta}) + \mu_1 \cdot (\gamma - \overline{\gamma}) + \mu_2 \cdot (\underline{\gamma} - \gamma)
\end{equation}
The Karush-Kuhn-Tucker (KKT) conditions for the optimal point $\gamma^*$ are:
\subsubsection{First order necessary conditions}
\begin{equation}
    \frac{\partial{\mathcal{L}}}{\partial{\gamma}} = -f_{\gamma}(\gamma^*; \Vec{\theta}) + \mu_1^* - \mu_2^* = 0
\end{equation}
\subsubsection{Complementary slackness conditions}
\begin{equation}
    \mu_1^* \cdot (\gamma^* - \overline{\gamma}) = 0 \; , \; \mu_2^* \cdot (\underline{\gamma} - \gamma^*) = 0
\end{equation}
\subsubsection{Dual problem feasibility}
\begin{equation}
    \mu_1^* \; , \mu_2^* \geq 0
\end{equation}
If one of the constraints $g^i$ is active at the nominal optimum, then its corresponding multiplier $\mu_i > 0$ by strict complementarity, and if inactive, then $\mu_i = 0$. In addition to the above necessary conditions, the following Strong $2^{nd}$ Order Sufficient Condition (SSC) guarantees that $\gamma^*$ is a local minimum of the NLP, even if it is non-convex. The Hessian of the Lagrangian must be positive definite on the null space of the Jacobian of active constraints $g^a$ \cite{buskens2001sensitivity}:
\begin{equation}
    \nu^\intercal \; \mathcal{L}_{\gamma \gamma} \; \nu > 0 \quad \forall \; \nu \; \neq \; 0 \quad \textrm{s.t.} \quad g^1_{\gamma}\; \nu = 0 \quad \textrm{or} \quad g^2_{\gamma}\;\nu = 0
\end{equation}
Since at most one constraint can be active, this implies that either $\nu = 0$ or $-\nu = 0$ and thus the SSC holds automatically if either the lower (\cref{eq:opt_prob1}) or upper (\cref{eq:opt_prob2}) bound is active. If neither constraint is active, then the SSC requires positive definiteness of $\mathcal{L}_{\gamma \gamma} = -f_{\gamma \gamma}$ over all possible values of $\gamma$ in the domain. It can be shown that the objective function $f(\gamma; \Vec{\theta})$ is concave in $\gamma$ as long as a specific condition holds on the parameters, price and travel times. For instance, if $R = \overline{u}$, this condition turns out to be:
\begin{align}
\label{eq:ssc}
    e^{-\lambda(\overline{u} - u_0)^{\beta}}\left(e^{-e^{-\lambda (\overline{u} - \underline{u})^{\beta} (-ln (p))^{\alpha}}} + \gamma \lambda \beta b_{sm}(\overline{u} - u_0)^{\beta -1}\right) \nonumber \\
    \leq - \left(e^{-e^{-\lambda (\overline{u} - \underline{u})^{\beta} (-ln (p))^{\alpha}}}\right)^2
\end{align}
This implies $\mathcal{L}_{\gamma \gamma} \geq 0$ and the NLP reduces to a convex optimization problem for which KKT conditions are sufficient for both local and global minima.
\subsection{\label{sec:local} Local sensitivity analysis}
Suppose the \textit{nominal} problem with assumed parameters $\theta_0$ has optimal tariff $\gamma^*_0$ and optimal objective value $f(\gamma^*_0; \theta_0) = f^*_0$. If the actual parameters turn out to be $\hat{\theta}$, we consider how this impacts the optimal tariff and optimal objective for the perturbed problem. Here, we consider cases where only one of the parameters is perturbed at a time while keeping the others fixed. This doesn't account for how interactions between parameters may influence the objective function. Local sensitivity analysis considers relatively small perturbations or uncertainties in the parameters for which the active set remains constant. Following \cite{buskens2001sensitivity}, local sensitivity differentials can be derived analytically in the neighbourhood of the nominal optimum operating point ($\gamma^*_0$,$\Vec{\theta_0}$), considering variations in a single parameter $\theta$:
\begin{gather}
    \begin{bmatrix} \frac{d\gamma^*}{d\theta} \\ \frac{d\mu^a}{d\theta} \end{bmatrix} = -
    \begin{bmatrix} \mathcal{L}_{\gamma \gamma}& {g^a}^\intercal \\ g^a & 0 \end{bmatrix}^{-1} \begin{bmatrix} \mathcal{L}_{\gamma \theta} \\ g^a \end{bmatrix}
    \label{eq:local_sens}
\end{gather}
where all the quantities are evaluated at the nominal values $\gamma^*_0$, $\theta_0$ and $\mu^a_0$. This gives us local sensitivity derivatives of both the optimal solution $\gamma^*(\theta)$ and multipliers $\mu^a$ corresponding to inequality constraints $g^a$ active at the nominal optimum. Furthermore, $\frac{d\mu^{ina}}{d\theta} = 0$ for all inactive constraint multipliers. If neither $g^1$ nor $g^2$ is active, we get:
\begin{align}
    \frac{d\gamma^*}{d\theta} & = -\mathcal{L}_{\gamma \gamma}^{-1} \mathcal{L}_{\gamma \theta}
\end{align}
The $1^{st}$ order sensitivity of the optimal objective yields \cite{buskens2001sensitivity}:
\begin{align}
    \frac{df^*}{d\theta}(\gamma(\theta);\theta)|_{\theta = \theta_0} & = \mathcal{L}_{\theta}(\gamma^*_0,\mu_0^a,\theta_0)
\end{align}
\subsection{\label{sec:Taylor} Real-time approximations by Taylor expansions}
We can also approximate the perturbed optimal solution $\gamma^*(\theta)$ and objective function $f^*$ in the neighbourhood of the nominal optimum, using $1^{st}$ order Taylor expansions about this operating point:
\begin{align}
    \gamma^*(\theta) & = \gamma^*_0 + \frac{d\gamma^*}{d\theta}(\theta_0) (\theta - \theta_0) \\
    f^*(\theta) & = f^*_0 + \frac{df^*}{d\theta}(\theta_0) (\theta - \theta_0)
\end{align}
We can measure the quality of this approximation by comparing it with the $2^{nd}$ order expansion \cite{buskens2001sensitivity}:
\begin{align}
\label{eq:2ndorder}
    f^*(\theta) & = f^*_0 + \frac{df^*}{d\theta}(\theta_0) (\theta - \theta_0) + \frac{d^2f^*}{d\theta^2}(\theta_0) (\theta - \theta_0)^2 \nonumber \\
    \frac{d^2f^*}{d\theta^2}(\theta_0) & = \mathcal{L}_{\gamma \gamma}(\theta_0) \left(\frac{d\gamma^*}{d\theta}(\theta_0)\right)^2 + 2\mathcal{L}_{\gamma \theta}\frac{d\gamma^*}{d\theta}(\theta_0) + \mathcal{L}_{\theta \theta}
\end{align}
\subsection{\label{sec:domain} Prediction of local domain}
The analysis in \cref{sec:local} assumes that the perturbation does not alter the set of active constraints. We can estimate the largest allowable magnitude of such parameter changes that still preserves the active set. $1^{st}$ order estimates for the Lagrange multipliers are used to determine when their corresponding constraints become either slack or tight, as in \cref{eq:domain1} and \cref{eq:domain2} respectively.  
\begin{enumerate}
    \item When a constraint leaves the active set, its non-zero multiplier becomes zero: 
    \begin{align}
    \label{eq:domain1}
        \mu^a(\theta) & \approx \mu^a(\theta_0) + \frac{d\mu^a}{d\theta}(\theta_0)(\theta - \theta_0) = 0 \nonumber \\
        (\hat{\theta} - \theta_0)_{max} & = \Delta \theta_{max} = - \frac{\mu^a(\theta_0)}{\frac{d\mu^a}{d\theta}(\theta_0)}
    \end{align}
    \item When a constraint enters the active set, it becomes tight and equal to zero.
    \begin{align}
    \label{eq:domain2}
        g_{ina}(\gamma;\theta) & \approx g_{ina}(\gamma^*_0;\theta_0) + \frac{dg_{ina}}{d\theta}(\gamma^*_0,\theta_0)(\theta - \theta_0) = 0 \nonumber \\
        (\hat{\theta} - \theta_0)_{max} & = \Delta \theta_{max} = - \frac{g_{ina}(\gamma_0;\theta_0)}{\frac{dg_{ina}}{d\theta}(\gamma_0;\theta_0)}
    \end{align}
\end{enumerate}
The size of this computed local domain represents the maximum allowable relative perturbation that does not change the active constraint set, thus allowing us to determine whether the analytical solution will be accurate.

\subsection{\label{sec:global} Global sensitivity analysis}
Global methods study the effects of varying multiple parameters simultaneously and relatively larger perturbations that cause the active set to change. In our CPT behavioral model, this could mean large errors in parameters or fundamentally misclassifying users, for example, assuming a given passenger to be loss averse ($\lambda > 1$) when in fact they are not ($\lambda < 1$). It is generally not possible to obtain explicit sensitivity derivatives at such points. However, an iterative scheme can be used to calculate directional derivatives for the optimal tariff and value function \cite{buskens2001sensitivity}:
\begin{enumerate}
    \item Calculate the initial optimum ($\gamma^*, \mu^a$) and sensitivity differentials $\left(\frac{d\gamma^*}{d\theta},\frac{d\mu^a}{d\theta}\right)$ at the nominal value $\theta_0 = \theta_0^1$.
    \item Compute the local domain as in \cref{sec:domain} and the perturbed parameter $\theta_0^2$ that disturbs the active set.
    \item Calculate sensitivity differentials at $\theta_0^2$ and update the active set to calculate $1^{st}$ order changes:
    \begin{align}
        \Delta \gamma^*& = \frac{d\gamma^*}{d\theta}(\theta_0^1)(\theta_0^2 - \theta_0^1) +  \frac{d\gamma^*}{d\theta}(\theta_0^2)(\theta - \theta_0^2) \\
        \Delta \mu^a & = \frac{d\mu^a}{d\theta}(\theta_0^1)(\theta_0^2 - \theta_0^1) +  \frac{d\mu^a}{d\theta}(\theta_0^2)(\theta - \theta_0^2)
    \end{align}
    \item Compute new optimal solutions and multipliers, as well as $1^{st}$ and $2^{nd}$ order approximations of 
    $f^*$.
    \item Repeat steps (1)-(4) whenever the active set is detected to change with incrementally larger perturbations.
\end{enumerate}
Thus, while local sensitivity analysis constructs a linear approximation around the nominal operating point, the global sensitivity analysis creates a \textit{piecewise} linear approximation with discontinuities at points where the active set changes. 
\subsection{Mismatch loss}
Parameter estimation errors cause losses in the objective function, resulting from a mismatch between parameters assumed while designing the dynamic price ($\widetilde{\theta}$) versus the true but unknown, behavioral model of passengers ($\theta_{true}$).
\begin{align}
    \Delta f & = f(\gamma^*_{true}; \theta_{true}) - f(\widetilde{\gamma}^*;\theta_{true}) \label{eq:mismatch} \\
    \gamma^*_{true} & = \argmax_\gamma \; f(\gamma; \theta_{true}) \\
    \widetilde{\gamma}^* & = \gamma^*(\widetilde{\theta}) = \argmax_\gamma \; f(\gamma; \widetilde{\theta}) 
\end{align}
which implies that $f(\gamma^*_{true}; \theta_{true}) \geq f(\widetilde{\gamma}^*; \theta_{true})$.
\subsection{Numerical simulations}
Several assumptions need to be made to obtain reasonably accurate analytical solutions, for e.g, regarding the largest magnitude of allowed perturbations for a localized analysis to be valid, accounting for possible changes in the active set (\cref{sec:domain}), checking curvature and convexity (\cref{eq:ssc}) etc. These can be quite restrictive especially if we wish to consider larger perturbations and uncertainties in parameters. However, when the above mentioned assumptions do hold, analytical methods can be advantageous and much faster since most of the calculations can be done offline. 

In addition, numerical approaches can also be used. The updated optimal solutions, value functions and mismatch losses can be computed using solvers in MATLAB's Global Optimization Toolbox, such as \texttt{Global Search} - which provides fast, proven quadratic convergence to local optima for such smooth problems using gradient-based methods. Results can be obtained by artificially constructing sensible travel scenarios (for both the SMoDS and the alternative) and repeatedly solving the NLP in eq. (\ref{eq:exp_rev})-(\ref{eq:opt_prob2}) under both nominal and perturbed conditions. Numerical results can be applied more generally and provide a benchmark against which we can measure the accuracy of analytical approximations. On the other hand, such simulation-based methods are much more computationally expensive. This brute-force method is feasible here since the problem size and associated computational burden are relatively small. However, it may not be practical for larger, higher dimensional problems with more constraints. As the number of passengers and trips increases, both the problem size and dimensionality increase since there's a dynamic tariff associated with each trip (or travel scenario) and each passenger has a unique set of user-specific CPT parameters describing their behavior. Furthermore, such numerical solvers generally do not provide proven, theoretical guarantees for convergence to global minima. However, this can be achieved practically by modifying the solver settings as needed. By rerunning the optimization using a large number of local solvers, \texttt{Global Search} can then return global minimizers. In our simulations, several termination criteria of these solvers were tweaked to reach (nearly) globally optimally solutions. These settings included stopping tolerances on gradients, step sizes, objective function values and finite difference settings etc.

\section{RESULTS AND DISCUSSION}
\subsection{Key insights from representative scenarios}
\begin{figure}[H]
  \centering
  \begin{subfigure}[b]{0.49\linewidth}
    \includegraphics[width=\linewidth]{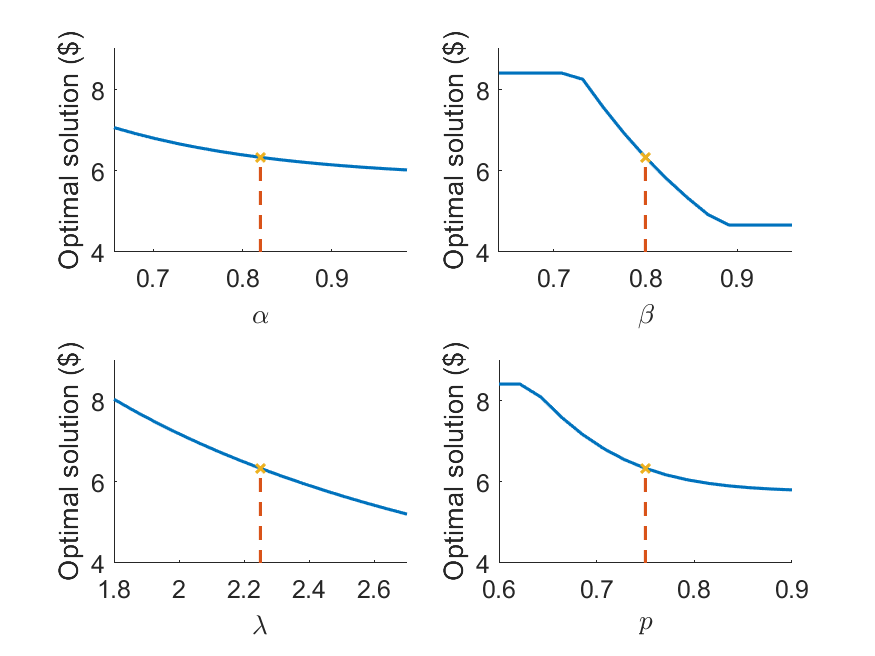}
    \caption{Variation in optimal solution $\gamma^*$ (SMoDS tariff).}
  \end{subfigure}
  \begin{subfigure}[b]{0.49\linewidth}
    \includegraphics[width=\linewidth]{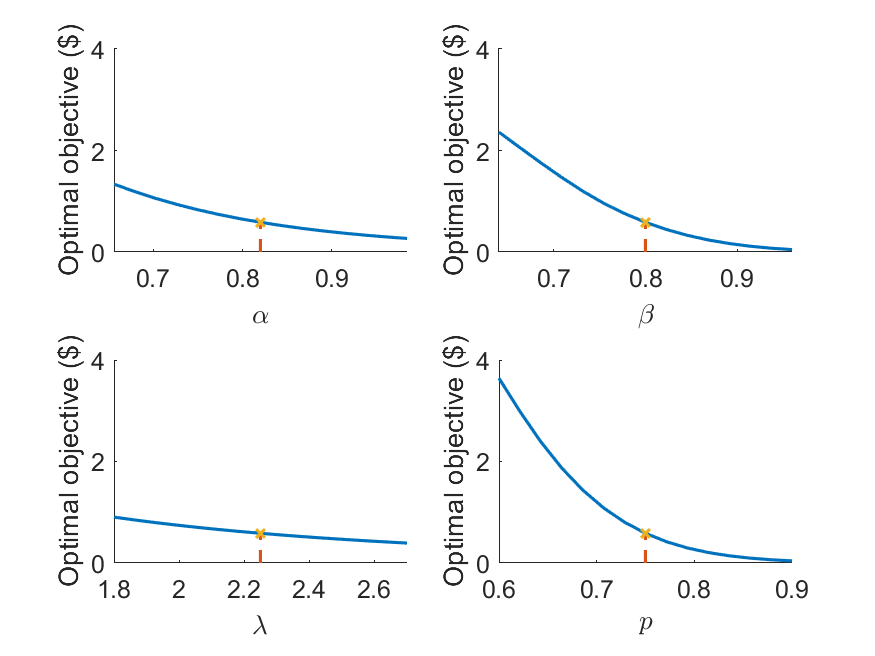}
    \caption{Variation in optimal objective $f^*$ (maximum expected revenue).}
  \end{subfigure}
  \begin{subfigure}[b]{0.49\linewidth}
    \includegraphics[width=\linewidth]{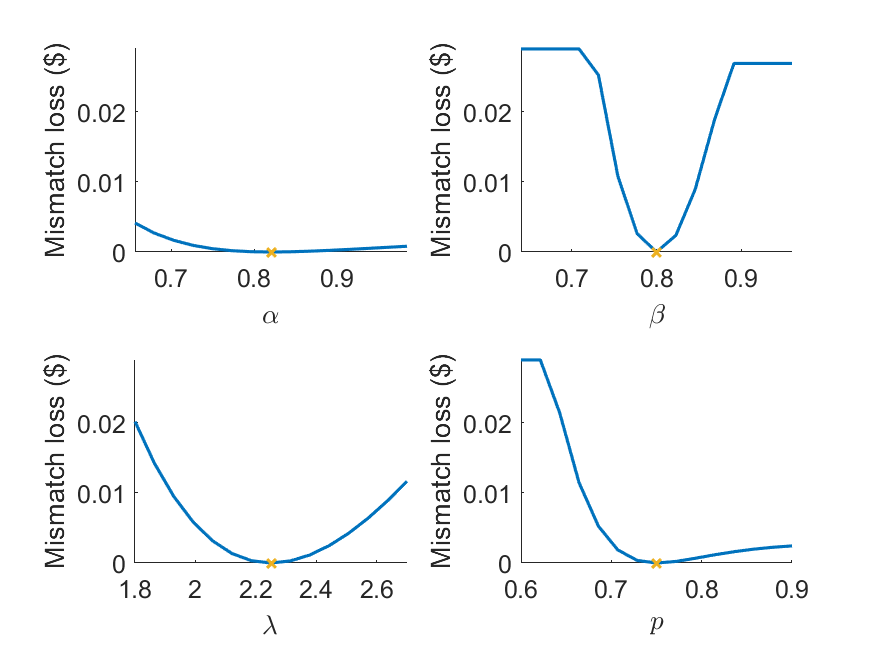}
    \caption{Variation in mismatch loss $\Delta f$.}
  \end{subfigure}
  \caption{\label{fig:S1} Sensitivity results for $S1$.}
\end{figure}

Both analytical and numerical results were obtained for more than 100 randomly generated travel scenarios. These were created by varying $u_0$, $\underline{x}$, $\overline{x}$, $b_{sm}$, $\underline{\gamma}$, $\overline{\gamma}$ while ensuring that the resulting choice set was valid and involved fair comparisons between the uncertain SMoDS and the certain alternative travel option, i.e., $\underline{u} \leq u_o \leq \overline{u}$. The objective utilities of the SMoDS and alternative can take positive or negative values according to \cref{eq:utilities}. The sensitivity analysis was performed for each of these scenarios using nominal parameter values $\alpha_0 = 0.82$, $\beta_0 = 0.8$ and $\lambda_0 = 2.25$ estimated in \cite{guan2019cumulative}. The probability of the worst-case SMoDS outcomes was set as $p_0 = 0.75$. In order to compute the mismatch loss, the passenger's true behavioral parameters were set equal to their initial values at the nominal optimum, i.e., $\theta_{true} = \theta_0$ in \cref{eq:mismatch}. All five of these model parameters are then varied by as much as $\pm 20\%$. Only five select scenarios are presented here to show some distinct behaviors and trends, as summarized in \cref{tab:scenarios}. For all scenarios, the best-case SMoDS outcome was used as reference, i.e., $R = \overline{u}$. In this case, both the alternative $A$ and the worst case SMoDS outcome are perceived as losses by the passenger. The acceptance probability can then be derived using equations (\ref{eq:value_func})-(\ref{eq:subj_prob}):
\begin{align}
    A_R^s & = V(u_0) = -\lambda (\overline{u} - u_0)^{\beta} \\
    U_R^s & = \pi(p) V(\underline{u}) = -\lambda e^{-(-ln(p))^\alpha} (\overline{u} - \underline{u})^{\beta} \\
    f =  \gamma p_R^s & = \frac{\gamma}{1 + e^{\lambda\left(e^{-(-ln(p))^\alpha} (\overline{x} - \underline{x})^\beta - (\overline{x} + b\gamma - u_0)^\beta\right)}} \label{eq:interpret}
\end{align}
\begin{table}[H]
\caption{Selected representative scenarios.}
\label{tab:scenarios}
\centering
\begin{tabular}{@{}ccccccc@{}}
\toprule
\textbf{Scenario} & $u_0 [-]$ & $b_{sm} [\$^{-1}]$ & $\underline{\gamma} [\$]$ & $\overline{\gamma} [\$]$ & $\underline{x} [-]$ & $\overline{x} [-]$ \\ \midrule
$S1$   &  8.17       &     -0.14     &   4.66   &     8.41  &   2.46 &  15.45                            \\
$S2$   & -7.62   &  -0.08  & 2.04  & 19.26 & -8.67  &  -5.15 \\
$S3$  & -2.54    &    -0.72  &  4.12   &   12.99   &               0.32     &     10.98  \\
$S4$   &    9.51  &   -0.40   &    1.11  &  13.49  &  1.06 & 24.36   \\
$S5$   &   9.55  &   -0.04    &   4.24  &   7.92  &  4.41  &  12.92  \\
\bottomrule
\end{tabular}
\end{table}
\begin{table}[H]
\caption{Sensitivity differentials of the optimal solution and value function near the nominal operating point, all evaluated at their respective nominal parameter values.}
\label{tab:differentials}
\centering
\begin{tabular}{@{}cccccccccccc@{}}
\toprule
\textbf{Scenario} & $\frac{d\gamma^*}{d\alpha}$ & $\frac{d\gamma^*}{d\beta}$ & $\frac{d\gamma^*}{d\lambda}$ & $\frac{d\gamma^*}{dp}$ & $\frac{df^*}{d\alpha}$ &  $\frac{df^*}{d\beta}$ & $\frac{df^*}{d\lambda}$ & $\frac{df^*}{dp}$\\ \midrule
$S1$  & -2.8 & -24.1 & -2.9 & -8.7 &  2.9    &    6.9    &  0.5  &  8.9  \\
$S2$  & -3.9 & -20.8 & -4.2 & -12.1 & 4.3  &   8.4   &    0.9    &   12.9           \\
$S3$ & -4.5 & -3.1 & 0.2 & -13.8 & 4.6 &  -0.8 &  -0.7  &   14.1                                              \\
$S4$   & -4.3 & -15.7 & -1.3 & -13.3 &  4.4   &  6.9  & 0.3  &  13.5 \\
$S5$   & -27.7 & -99.7 & -12.2 & -84.5 & 1.2  &  3.8   &  0.4  & 3.6 \\ 
\bottomrule
\end{tabular}
\end{table}
\begin{table}[H]
\caption{Valid domain for local sensitivity analysis (in \%).}
\label{tab:domains}
\centering
\begin{tabular}{@{}ccccc@{}}
\toprule
\textbf{Scenario} & $|\Delta \alpha_{max}|$ & $|\Delta \beta_{max}|$ & $|\Delta \lambda_{max}|$ & $|\Delta p_{max}|$ \\ \midrule
$S1$  &   71.7    &   8.7     &  24.6  &  25.7       \\
$S2$ &  250.8    &   48.9  &    85.6   &  89.9    \\
$S3$ &  88.8  &  131.9   &  706.9  &   31.8         \\
$S4$  & 51.9   &  14.8  &  63.5 & 18.6  \\
$S5$  & $3\times 10^{-6}$  & $10^{-6}$   &  $3\times 10^{-6}$  & $10^{-6}$ \\
\bottomrule
\end{tabular}
\end{table}
The plots in \cref{fig:S1} display variations in the optimal tariff, objective function value (maximum expected revenue) and mismatch loss for scenario $S1$. The variations in optimal solution and mismatch loss sometimes become flat w.r.t. to certain parameters as seen with $\beta$ and $p$ in \cref{fig:S1}. This happens because $\gamma^*$ hits either the lower or upper bound and doesn't change further unless the active set changes yet again due to even larger perturbations. We now present three major implications from analyzing these scenarios. 
\begin{figure}[H]
  \centering
  \begin{subfigure}[b]{0.49\linewidth}
    \includegraphics[width=\linewidth]{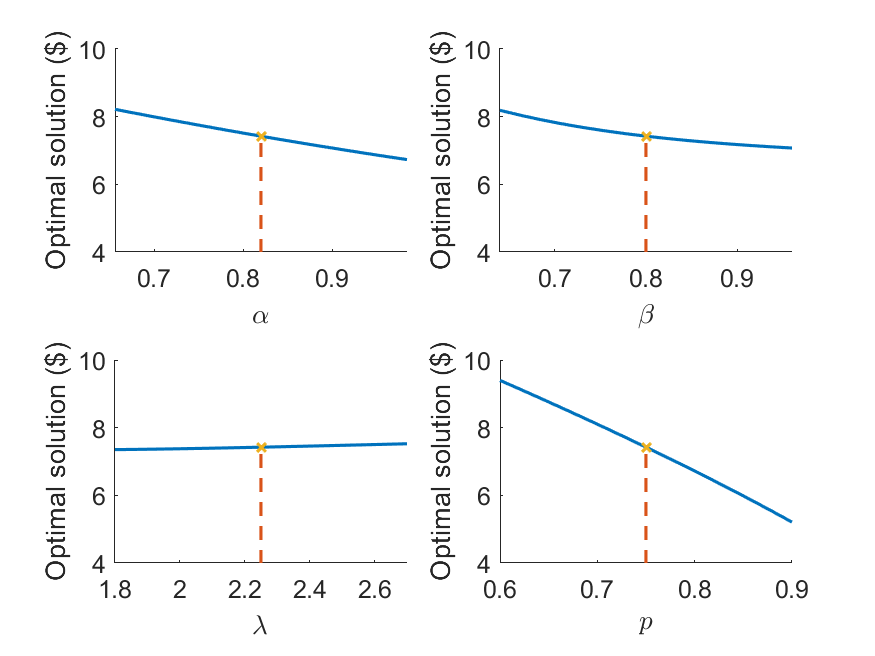}
    \caption{Variation in optimal solution.}
  \end{subfigure}
  \begin{subfigure}[b]{0.49\linewidth}
    \includegraphics[width=\linewidth]{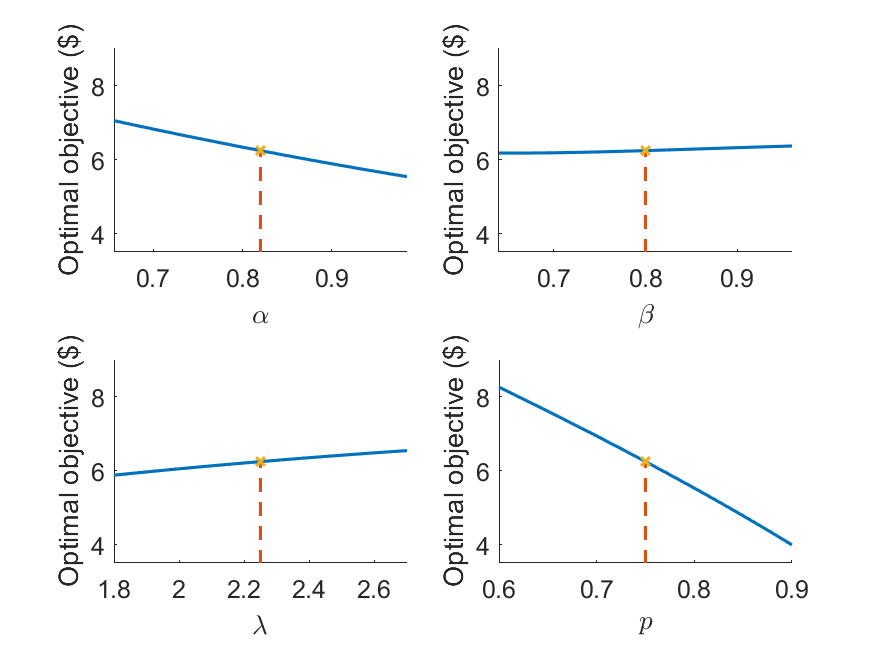}
    \caption{Variation in optimal objective.}
  \end{subfigure}
  \caption{\label{fig:S3} Sensitivity results for $S3$.}
\end{figure}
\subsubsection{\label{sec:general} General sensitivity trends for different parameters}
\begin{figure}[H]
  \centering
  \begin{subfigure}[b]{0.49\linewidth}
    \includegraphics[width=\linewidth]{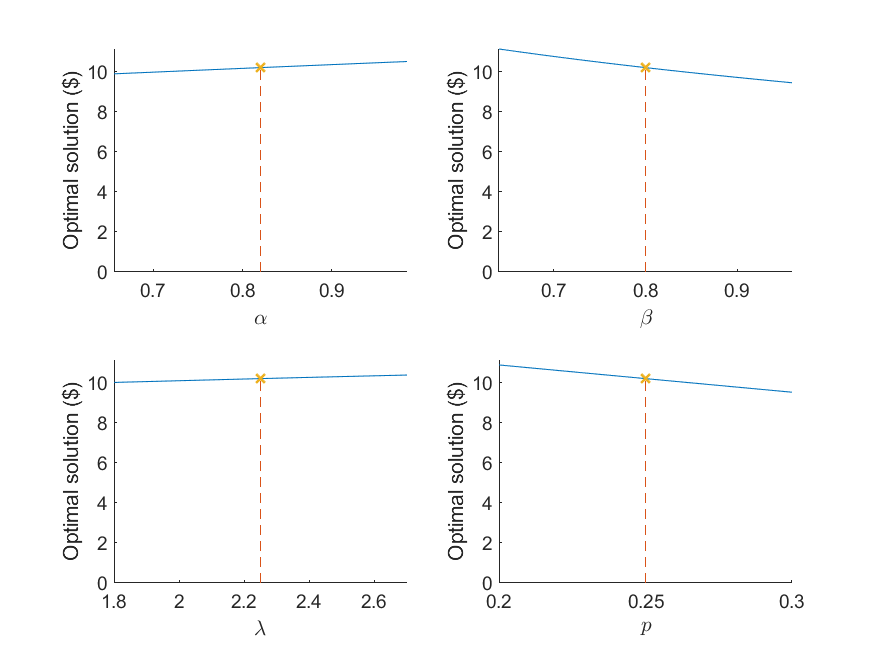}
    \caption{Variation in optimal solution (SMoDS tariff).}
  \end{subfigure}
  \begin{subfigure}[b]{0.49\linewidth}
    \includegraphics[width=\linewidth]{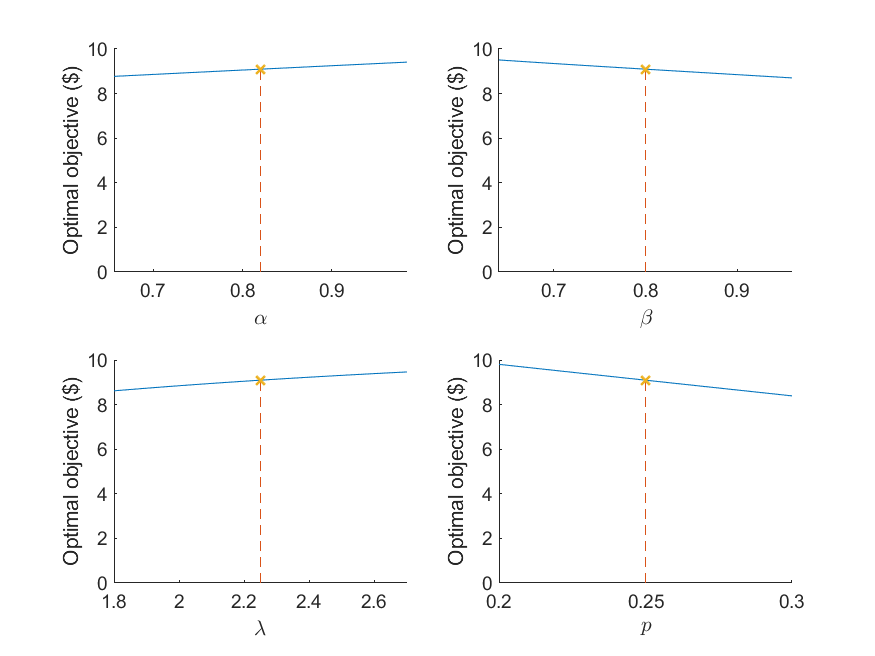}
    \caption{Variation in optimal objective.}
  \end{subfigure}
  \caption{\label{fig:S6} Sensitivity results for $S6$ with $p = 0.25$.}
\end{figure}

Across all the scenarios studied, the optimal solution was found to be monotonically decreasing with respect to both $\alpha$ and $p$. This makes intuitive sense since as $p$ increases, the likelihood of the best-case outcome decreases while the worst-case outcome becomes more likely. This makes the SMoDS relatively less attractive compared to the certain alternative $A$ and thus the optimal SMoDS tariff $\gamma^*$ must be reduced accordingly for it to remain competitive. Similarly, an increase in $\alpha$ indicates weaker probability distortion (i.e. under weighting high probability events and over weighting low probability events). Thus, as $\alpha$ increases, the user underestimates the likelihood of the more likely, worst-case outcome occurring to a lesser extent or equivalently overestimates the rarer, best-case outcome less. This reduces the SMoDS' relative attractiveness as well as $\gamma^*$. The opposite (monotonically increasing) trend would be observed with respect to $\alpha$ if best-case SMoDS outcome occurred with higher probability instead of the worst-case outcome, i.e., if $p < 0.5$. This was verified by also testing a few scenarios with $p_0 = 0.25$ instead of $p_0 = 0.75$, one of which is shown in \cref{fig:S6}.
The variations in optimal objective are less straightforward to predict. Since $f^*(\gamma^*; \Vec{\theta}) = \gamma^* \cdot p_R^s(\gamma^*; \Vec{\theta})$ and $p_R^s$ is strictly monotonically decreasing in $\gamma$ \cite{guan2019cumulative}, the variation in expected revenue due to parameter changes will depend on how much the acceptance probability rises in response to a fall in $\gamma^*$ and vice versa.

\subsubsection{High degree of scenario specificity}
Although we can make some broad statements about trends w.r.t. $\alpha$ and $p$, the same cannot be said for other parameters. Changes in $\alpha$ and $p$ only affect the subjective utility of the SMoDS but variations in $\beta$ and $\lambda$ affect both the SMoDS as well as the alternative (\cref{eq:interpret}). Thus, either a monotonically increasing or decreasing trend could be obtained depending on the relative magnitudes of $e^{-(-ln(p))^\alpha} (\overline{x} - \underline{x})^\beta$ and  $(\overline{x} + b\gamma^* - u_0)^\beta$ in \cref{eq:interpret}. In other words, the behavior of the optimal solution in response to parameter perturbations depends on (1) the spread between the two possible SMoDS outcomes and (2) how this compares with the alternative. For instance, if $e^{-(-ln(p))^\alpha} (\overline{x} - \underline{x})^\beta > (\overline{x} + b\gamma^* - u_0)^\beta \implies U_R^s < A_R^s < 0$, i.e., the passenger perceives the worst-case SMoDS outcome as a greater loss (lower subjective utility) than the alternative. An increase in $\lambda$ would make the user more averse to losses and thus persuade them away from the SMoDS in favor of the alternative. This reduces the relative attractiveness of SMoDS, causing $\gamma^*$ to fall. The opposite effect occurs if $e^{-(-ln(p))^\alpha} (\overline{x} - \underline{x})^\beta < (\overline{x} + b\gamma^* - u_0)^\beta$ and the optimal tariff $\gamma^*$ increases monotonically with $\lambda$, as seen for scenario $S3$ in \cref{fig:S3}.
\begin{figure}[H]
  \centering
  \begin{subfigure}[b]{0.49\linewidth}
    \includegraphics[width=\linewidth]{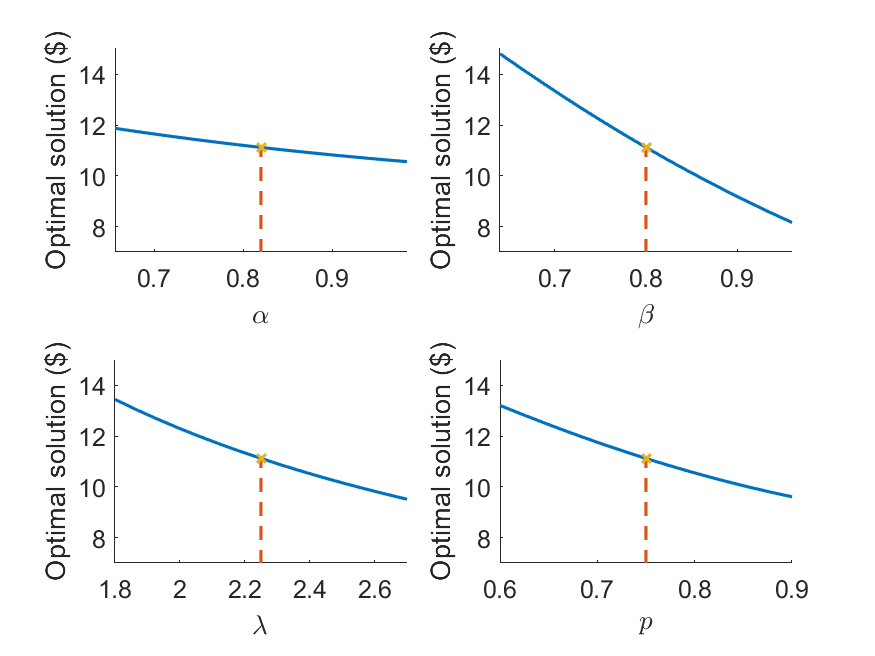}
    \caption{Variation in optimal solution.}
  \end{subfigure}
  \begin{subfigure}[b]{0.49\linewidth}
    \includegraphics[width=\linewidth]{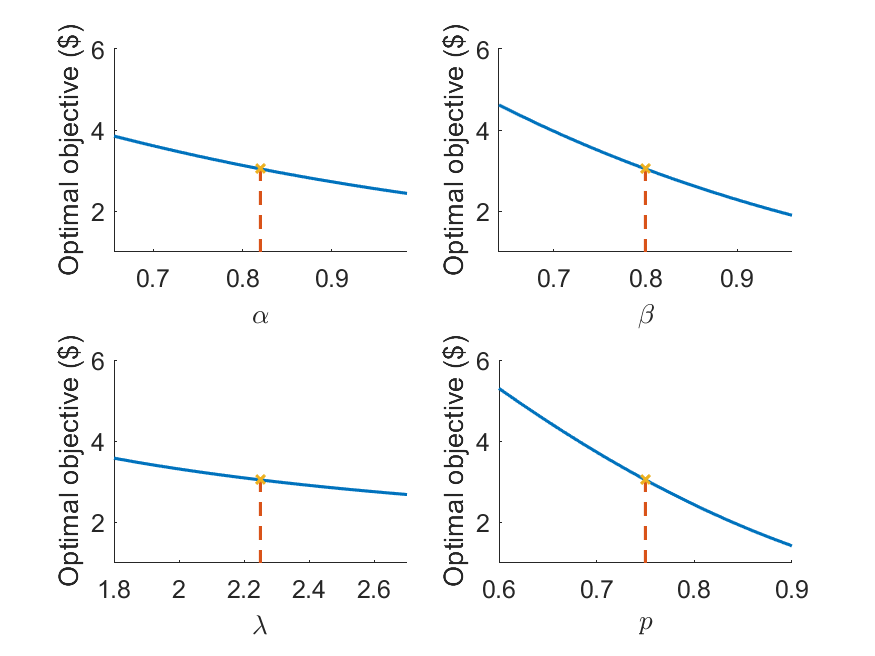}
    \caption{Variation in optimal objective.}
  \end{subfigure}
  \begin{subfigure}[b]{0.49\linewidth}
    \includegraphics[width=\linewidth]{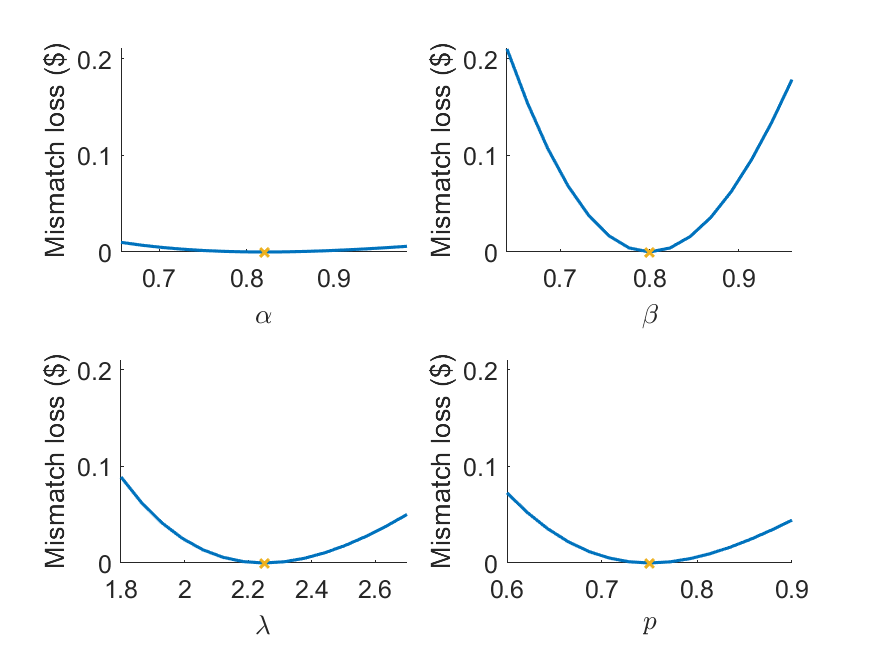}
    \caption{Variation in mismatch loss.}
  \end{subfigure}
  \caption{\label{fig:S2} Sensitivity results for $S2$.}
\end{figure}
Similarly, across a vast majority of the scenarios studied, $\gamma^*$ was found to be decreasing with $\beta$. This agrees with intuition since $0 < \beta < 1$ describes the rate at which sensitivity (to both losses and gains) decreases as the outcome moves farther away from the reference. Thus, as $\beta$ increases, the passenger becomes relatively more sensitive to losses even away from $R = \overline{u}$. Since $\underline{u} \leq u_0 \leq \overline{u}$, the passenger is now hurt more by the worst-case SMoDS outcome than before because the diminishing sensitivity effect is less influential. Thus, the optimal tariff $\gamma^*$ decreases to compensate for this.  

Furthermore, although we can sometimes predict the general direction of the variation as in \cref{sec:general}, the exact nature of the trend observed, i.e., the magnitude and rate of resulting changes, depends significantly upon the specific choice set and travel scenario under consideration. These are determined by the objective utility of the certain alternative ($u_0$), travel time outcomes possible with the SMoDS ($\underline{x}$, $\overline{x}$), disutility associated with the SMoDS tariff ($b_{sm}$) and the bounds placed on the dynamic price ([$\underline{\gamma},\overline{\gamma}$]). Even after experimenting with a large number of randomized scenarios, we were unable to obtain very specific scenario-agnostic insights that could be generalized. This is likely a consequence of the strongly nonlinear nature of the objective function. Thus, we concluded that the exact nature of variations in the optimal solution, value function and mismatch loss are determined to a large extent by the particular characteristics of each family of scenarios. These include various features like trip duration and distance covered, the passenger's utility functions (value of time) for each mode option and each leg of the trip, nominal tariff level, etc. the In addition, the trends observed are likely to be affected by the reference level ($R$) used as well as the nominal values of the CPT parameters. For instance, we obtain almost linear variations for travel scenario $S2$ in \cref{fig:S2}, which is quite different from the behavior seen for $S1$ in \cref{fig:S1}. The relative sensitivities w.r.t. to each parameter can also vary quite significantly depending on the specific travel scenario as well as the nominal parameters, as indicated by the sensitivity differentials in \cref{tab:differentials}. 

\subsubsection{Robustness under certain scenarios}
The local sensitivity domain was found to be quite large for most combinations of nominal parameters and travel scenarios, as shown by $S1$-$S4$ in \cref{tab:domains}. This implies that a relatively large perturbation in parameters is needed to alter the nominal active set. Such a property can be exploited to set more robust dynamic tariffs. If the ride offer and travel scenarios can be designed such that the optimal tariff is at either the upper ($\overline{\gamma}$) or lower ($\underline{\gamma}$) bound, then $\gamma^*$ would remain unchanged even with large errors or misconceptions in the parameter values assumed for the passenger. Although we assumed the tariff bounds to be exogenous to the travel scenarios in this study, the choice of this set $[\underline{\gamma},\overline{\gamma}]$ will likely involve another optimization problem in itself. In practical applications, the bounds specified on the tariff would be determined depending on the given scenario and would need to be updated for each new passenger and trip. This aspect will be looked into further as part of future work. 

\subsection{\label{sec:num_ana_comp} Comparison of numerical results and analytical solutions}
\begin{figure}[H]
  \centering
  \begin{subfigure}[b]{0.49\linewidth}
    \includegraphics[width=\linewidth]{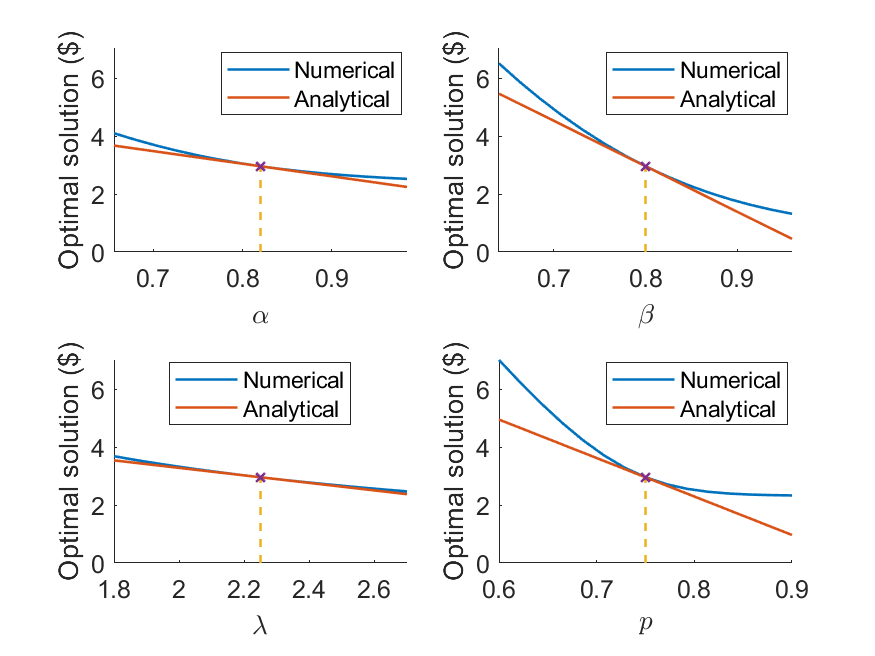}
    \caption{\label{fig:SolComparison1} Optimal price variation in $S4$.}
  \end{subfigure}
  \begin{subfigure}[b]{0.49\linewidth}
    \includegraphics[width=\linewidth]{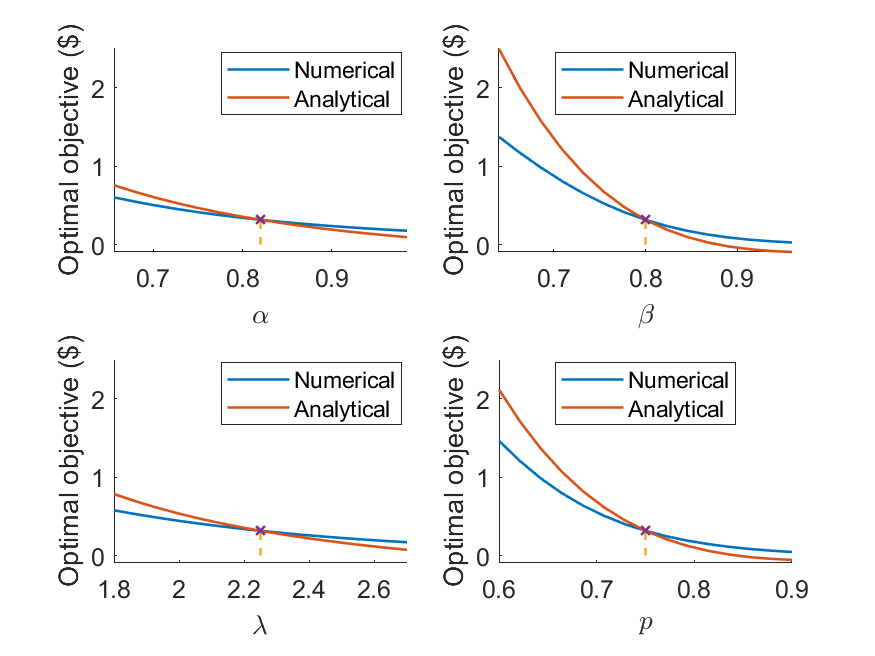}
    \caption{\label{fig:CostComparison12} Optimal objective variation in $S5$.}
  \end{subfigure}
  \caption{Analytical (local) vs numerical results.}
\end{figure}

In general, for a majority of the scenarios considered, linearized approximations using $1^{st}$ order Taylor expansions provide a reasonable estimate of the optimal solution for small \textit{local} perturbations (e.g. $\pm 20\%$) in the neighbourhood of the nominal parameter value. The analytical solutions obtained using a local sensitivity analysis, for both the optimal tariff and value function are very close to the exact values obtained numerically through repeated optimization, with $< 0.5 \%$ error for most scenarios.  However, larger discrepancies were obtained for a small subset of scenarios. For each of these scenarios, at least one of the following conditions was found to be true: 
\begin{enumerate}
    \item The variations in optimal tariff and optimal objective are clearly nonlinear and display significant curvature even for small perturbations from the nominal parameter values. In such cases, the error in the $1^{st}$ order linear approximation of the optimal solution is no longer negligible, as can be seen for parameter $\beta$ and $p$ in \cref{fig:SolComparison1}. However, for most of the scenarios, we observe close to linear variations in optimal objective for small to moderate perturbations. This makes intuitive sense since all four parameters ($\alpha$, $\beta$, $\lambda$ and $p$) are of the order of 1. Thus, the magnitude of deviations $\Delta \theta = \theta - \theta_0 < 1$ for small to medium perturbations, implying that the $2^{nd}$ and higher order terms in \cref{eq:2ndorder} can usually be ignored for scenarios where curvature is not important i.e. $\frac{d^2f^*}{d\theta^2}(\theta_0) << \frac{df^*}{d\theta} (\theta_0)$. \\
    
    \item  The local domain ($\Delta \theta_{max}$) is very small as for scenario $S5$ in \cref{tab:domains}, implying that even a slight perturbation in the parameter would change the active constraint set. This makes the optimal solution more sensitive to parameter uncertainty since $\gamma^*$ either (1) switches between the two bounds, (2) hits one of the bounds from the interior of $[\underline{\gamma},\overline{\gamma}]$ or (3) leaves one of the bounds for the interior of the set. This causes stronger nonlinear behavior in the optimal value function too and makes a local analysis insufficient, as seen for $\beta$ and $p$ in \cref{fig:CostComparison12}. This also agrees with intuition since linearized analysis is valid only for small deviations from the operating point.
    \end{enumerate}

For such cases, a global analytical method like that described in \cref{sec:global} could be used instead to obtain more accurate analytical solutions through better, piecewise linear approximations for the optimal tariff. This will be considered as part of future work.

\section{CONCLUDING REMARKS}
\subsection{Conclusions}
In this paper, we developed an analytical and numerical framework to quantify and predict the sensitivity of CPT-based SMoDS dynamic pricing schemes, to uncertainties in the behavioral model parameters for individual users. A few general trends could be identified for some of the parameters that held across scenarios. These align with the key axioms regarding subjective decision making of human beings that are examined in CPT \cite{tversky1992advances}. However, it was found that the exact nature of variations in the optimal tariff, expected revenue and mismatch loss are strongly influenced by specific attributes of the SMoDS trip offer as well as that of the alternative travel option under consideration. Quantifying the sensitivity using metrics like the local domains and sensitivity differentials also provides a tool to rank and prioritize the parameters in terms of their influence on the desired objective function. The analytical solutions based on local sensitivity approaches provide sufficiently accurate results and agree with numerical simulations for most travel scenarios. However, a global sensitivity analysis is found to be more appropriate in some special cases. 

\subsection{Future Work and other applications}
One important area that needs to be studied in greater detail is the dependence of the sensitivity and robustness trends on the various features of the travel scenarios. Although we considered a large number of random scenarios, a large portion of the parameter space still remains to be explored. In reality, there are an infinite number of possible scenarios corresponding to different mode-choice models, passengers, travel alternatives and trip attributes. Future work will also look into extending these results to more general and complex cases with fewer assumptions. This could entail a larger number of travel alternatives, possible SMoDS outcomes and following other probability distributions for SMoDS travel times. While this paper focused on using the best-case outcome as a dynamic reference, it would be interesting to look at other reference types as well as objective functions other than expected revenue. 

This study solely focused on single parameter perturbations (i.e. changing one factor at a time), but a similar approach can also be used while varying multiple parameters simultaneously. One could argue that the parameters $\alpha$, $\beta$ and $\gamma$ are uncorrelated and independent for a given individual since they represent distinct risk attitudes and behavioral patterns - this insight could be leveraged to simplify the analysis for multiple parameter perturbations. Another interesting area for exploration is scaling up this analysis to the population level. So far, the sensitivity has been studied only with respect to one passenger on a single trip. Simulating larger numbers of trips for an entire population of passengers in a city would give a better idea of the full impact of parameter uncertainties on overall revenue, ridership and performance of ridesharing fleets as a whole. In addition, future work will also focus on implementing methods for global sensitivity analysis, as mentioned earlier in \cref{sec:num_ana_comp}.

The framework and methods presented here can also be extended to applications beyond sensitivity analysis and managing parameter uncertainty. This paper focused on understanding the CPT model's sensitivity at a single point in time. However, in a realistic setting, these parameters fluctuate continuously. Thus, SMoDS pricing schemes would ideally involve closed-loop dynamic optimization for purposes of transactive control \cite{annaswamy2018transactive}. Fast sensitivity-based solution updates could be used to optimize tariffs in real-time while also efficiently handling uncertainty \cite{kadam2004sensitivity}. Such strategies would be able to not only maximize desired objective functions for the SMoDS system, but also learn the true passenger behavioural model parameters over time, potentially using their responses to each trip offer as feedback signals.

\subsection{Acknowledgements}
This work was supported by the Ford-MIT Alliance. The authors have no competing interests to declare.

\bibliographystyle{ieeetr}
\bibliography{AAmain}

\end{document}